\documentclass[12pt,british,a4paper]{amsart}

\usepackage{lmodern}
\usepackage[T1]{fontenc}
\usepackage{microtype}

\usepackage{babel,amssymb,amscd,url,stmaryrd,mathtools,fixltx2e}
\usepackage{paralist,csquotes,stackrel}
\usepackage[mathscr]{eucal}

\numberwithin{equation}{section}
\mathtoolsset{mathic=true}
\setdefaultenum{\upshape (i)}{}{}{}

\theoremstyle{plain}

\theoremstyle{definition}

\theoremstyle{remark}
\newtheorem{Example}[equation]{Example}

\newenvironment{example}{\begin{Example}\pushQED{\qee}}{\popQED\end{Example}}
\DeclareRobustCommand{\qee}{%
  \ifmmode \mathqee
  \else
    \leavevmode\unskip\penalty9999 \hbox{}\nobreak\hfill
    \quad\hbox{\qeesymbol}%
  \fi
}
\newcommand{\mathqee}{\quad\hbox{\qeesymbol}}
\newcommand{\qeesymbol}{\ensuremath\diamondsuit}

\newcommand{\Lie}[1]{\operatorname{\textsl{#1}}}
\newcommand{\lie}[1]{\operatorname{\mathfrak{#1}}}

\newcommand{\bmf}{\lie b}

\newcommand{\kf}{\lie k}
\newcommand{\h}{\lie h}
\newcommand{\n}{\lie n}

\newcommand{\tf}{\lie t}

\newcommand{\SL}{\Lie{SL}}
\newcommand{\SU}{\Lie{SU}}

\newcommand{\C}{{\mathbb C}}
\newcommand{\HH}{{\mathbb H}}

\newcommand{\R}{{\mathbb R}}
\newcommand{\Z}{{\mathbb Z}}

\newcommand{\symp}{{\sslash}} %extra brackets to force \mathord
\newcommand{\hkq}{{\sslash\mkern-6mu/}}

\DeclareMathOperator{\rank}{rank}

\usepackage{tikz} 

\usetikzlibrary{decorations.pathreplacing} 
\usetikzlibrary{decorations.pathmorphing} 
\usetikzlibrary{decorations.markings} 
\usetikzlibrary{arrows} 
\usetikzlibrary{shapes} 
\usetikzlibrary{matrix} 
\usetikzlibrary{positioning} 

\tikzset{cross/.style={cross out, draw=black}} 
\tikzset{circ/.style={circle,fill=white,draw=black}} 
\tikzset{hasse/.style={circle, fill,inner sep=2pt}} 
\tikzset{h/.style={circle, fill,inner sep=2pt}} 
\tikzset{ns/.style={circle, draw,inner sep=2pt}} 
\tikzset{dot/.style={circle,draw,fill=black}} 
\tikzset{gauge/.style={circle,draw}} 
\tikzset{flavor/.style={regular polygon,regular polygon sides=4, draw}} 
\tikzset{doublearrow/.style={ draw=black!75, color=black!75, thick, double distance=3pt, }} 
\tikzset{thirdline/.style={ draw=black!75, color=black!75, thick, }} 
\tikzset{middlearrow/.style={ decoration={markings, mark= at position 0.5 with {\arrow{#1}} , }, postaction={decorate} } } 
\tikzset{subquiver/.style={circle,draw,dashed,inner sep=0pt,minimum size=1cm}}

\begin{document}

\title{Symplectic duality and implosions}

\author[A.~Dancer]{Andrew Dancer}
\address[Dancer]{Jesus College\\
Oxford\\
OX1 3DW\\
United Kingdom} \email{dancer@maths.ox.ac.uk}

\author[A. Hanany]{Amihay Hanany}
\address[Hanany]{Theoretical Physics Group, The Blackett Laboratory,
Imperial College, Prince Consort Road, London SW7 2AZ, United Kingdom}
\email{a.hanany@imperial.ac.uk}

\author[F.~Kirwan]{Frances Kirwan}
\address[Kirwan]{New College\\
Oxford\\
OX1 3BN\\
United Kingdom} \email{kirwan@maths.ox.ac.uk}


\begin{abstract}
We discuss symplectic and hyperk\"ahler implosion and
present candidates for the symplectic duals of the universal hyperk\"ahler implosion
for various groups.
\end{abstract}

\maketitle

\section{Introduction}
Implosion is an abelianisation construction that originated in symplectic geometry \cite{GJS:implosion}.
and for which a hyperk\"ahler analogue was developed in a series of papers
\cite{DKS,DKS-Seshadri,DKS-twistor,
DKS-Arb}. In particular a complex-symplectic analogue of the
universal symplectic implosion for a compact simple
group was introduced, which in the $A_n$ case (ie the group $SU(n+1))$
is in fact hyperk\"ahler as a stratified space. The universal
implosion for $K$ carries a complex-symplectic
 action of $K_{\C} \times T_{\C}$ where
$T_{\C}$ is the complexification of the maximal torus $T$. In the $A_n$ case this is
the complexification of an action of $K \times T$ which preserves the
hyperk\"ahler structure (that is, it is isometric and triholomorphic).There is also
an action of $Sp(1)$ that rotates complex structures. 

This data suggests that there should be a symplectic dual of the implosion.
In this paper we present candidates for the symplectic duals in the
$A_n$ and $D_n$ cases, including some computational evidence.
We also include a discussion of implosions
and their links to quiver varieties and the Moore-Tachikawa category, which we hope will be of interest to
string theorists and algebraic geometers.

\subsubsection*{Acknowledgements.}

We thank BIRS for its hospitality during the workshop
``The analysis of gauge-theoretic moduli spaces'' in September 2017,
We thank Hiraku Nakajima for discussions during that workshop.

\section{Symplectic Implosion} \label{symplectic}

In this section
we review  the symplectic implosion construction of Guillemin, Jeffrey and
Sjamaar \cite{GJS:implosion}

The idea is that given a space $M$ with Hamiltonian 
action of a compact group $K$, one can form the imploded space $M_{\rm impl}$
with a Hamiltonian action of the maximal torus $T$ of $K$, such that
the symplectic reduction $M$ of $K$ agrees with the 
reduction of $M_{\rm impl}$ by $T$ as long as we reduce at levels in the
closed positive Weyl chamber. We can summarise this, using the usual
notation for symplectic quotients, as:
\[
M \symp_{\xi} K = M_{\rm impl} \symp_{\xi} T  \;\; : \;\; \xi \in
\bar{\tf}_{+}^*
\]
Fortunately the problem of constructing symplectic implosions can be
reduced to the case $M = T^*K$, which in this sense plays a universal role
for Hamiltonian spaces with $K$ action. The key point here is that
$T^*K$ has a Hamiltonian $K \times K$ action so when we form the implosion
$(T^*K)_{\rm impl}$ with respect to, say, the right $K$ action, the implosion
 has a $K \times T$ action, because the left $K$ action survives the
implosion process. Now
the implosion of a general symplectic manifold $X$ with Hamiltonian $K$-action
can be obtained by reducing $X \times (T^*K)_{\rm impl}$ by the diagonal $K$ action,
producing a space $X_{\rm impl}$ with $T$ action. The reduction of $X$ by $K$, at any element $\xi$ of a chosen closed positive Weyl chamber in the dual $\kf^*$ of the Lie algebra of $K$, coincides with the reduction of $X_{\rm impl}$ by $T$ at $\xi$.
In this sense the implosion abelianises the $K$ action on $X$.

The space $(T^*K)_{\rm impl}$ is referred to therefore as
the {\em universal symplectic implosion} for $K$. It is explicitly
constructed as a symplectic stratified space, by considering
the product $K \times \bar{\tf}_{+}^{*}$ of the group and the closed
positive Weyl chamber, and then performing certain collapsing 
operations as follows.

To motivate this, recall that the universal 
implosion should carry a Hamiltonian $K \times T$ action.
The reductions by $T$ at points in the closed positive Weyl
chamber should coincide with the reductions
of $T^*K$ by the right $K$ factor in the $K \times K$ action on $T^*K$.
These reductions are exactly the coadjoint orbits of $K$ : the $K$ action on these coadjoint orbits is induced by the left $K$ action on $T^*K$, or equivalently
by the $K$ factor in the $K \times T$ action on $(T^*K)_{\rm impl}$.

Now, for $K \times \bar{\tf}_{+}^{*}$
the $T$ moment map is projection onto the $\bar{\tf}_{+}^{*}$ factor
so the reduction at level $\xi$ is just $(K \times \{ \xi \})/T \cong K/T$.
This gives the correct picture for $\xi$ in the open Weyl chamber, but
not for $\xi$ in the lower-dimensional faces of the chamber.

If we stratify the product $K \times \bar{\tf}_{+}^{*}$ by the faces of the 
Weyl chamber, then the choice of stratum corresponds to a choice of
stabiliser $C$ for $\xi$, and the coadjoint orbit of $\xi$ is now $K/C$.
Therefore to obtain the coadjoint orbits on reduction by $T$, we must
 quotient each stratum by the commutator $[C,C]$. Now the
reduction by $T$ at level $\xi$ is 
$(K \times \{ \xi \})/T. [C,C] = K/{\rm Stab}_{K}(\xi)$ as required.

Hence the implosion is the symplectic stratified space obtained from
$K \times \bar{\tf}_{+}^{*}$ by stratifying by the faces of the Weyl chamber
and quotienting by the commutator of the stabiliser associated
to each stratum.
In particular no collapsing occurs on the open Weyl chamber
as $C$ is then abelian. This yields the top stratum $K \times \tf_{+}^{*}$.

\section{Nonreductive quotients} \label{nonreductive}
As often is the case with constructions in symplectic geometry, there
is an alternative description of the universal symplectic implosion in
terms of algebraic geometry.

We recall that geometric invariant theory (GIT) defines the quotient
$X/\!/G$ of an affine variety $X$ over $\C$ by the action of a complex
reductive group $G$ to be the affine variety
$\mathrm{Spec}(\mathcal{O}(X)^G)$ associated to the algebra
$\mathcal{O}(X)^G$ of $G$-invariant regular functions on $X$. This
is well-defined because in this situation the algebra $\mathcal{O}(X)^G$
is finitely generated. 

Moreover the inclusion of $\mathcal{O}(X)^G$ in $\mathcal{O}(X)$ induces a natural $G$-invariant morphism from $X$ to $X/\!/G$. When $G$ is reductive this morphism is always surjective, and points of $X$ become identified in $X/\!/G$ if and only if the closures of their $G$-orbits meet in the semistable locus of $X$. 

If $G$ is nonreductive then this picture can break down
because the algebra of invariants is not
necessarily finitely generated so $\mathrm{Spec}(\mathcal{O}(X)^G)$ need not
define an affine variety.
  Even if the algebra of
invariants is finitely generated, so that the GIT quotient
exists, the natural morphism $X \rightarrow X/\!/G$ is not necessarily 
surjective, and its image is in general not a subvariety of the GIT quotient but only a constructible subset \cite{DK} (ie a finite union of intersections
of open sets and closed sets).

It was shown in \cite{GJS:implosion} that the universal symplectic
implosion for a compact group $K$ can be identified with the
nonreductive GIT quotient $K_\C/\!/N$.  Here $K_{\C}$, the
complexification of $K$, is a complex affine variety, and $N$
denotes the maximal unipotent subgroup of $K_{\C}$.  Although $N$ is not reductive,
the algebra of invariants $\mathcal{O}(K_\C)^N$ is finitely generated
so $K_{\C} /\!/N$ exists as an affine variety. In fact $K_\C/\!/N$ 
may be viewed as the canonical
affine completion of the quasi-affine variety $K_\C/N$, which embeds
naturally as an open subset of $K_\C/\!/N$ with complement of
codimension at least two. The restriction map from
$\mathcal{O}(K_\C/\!/N)$ to $\mathcal{O}(K_\C/N)$ is thus an
isomorphism, and both algebras can be identified with the algebra of
$N$-invariant regular functions on $K_\C$. 

Moreover, there is a natural description of $K_{\C}/\!/N$ as a stratified space, where the strata may be identified
with $K_{\C}/[P,P]$ and $P$ ranges over the $2^{{\rm rank} K}$ standard
parabolics of $K_{\C}$. The top stratum, corresponding to choosing $P$ to
be the Borel subgroup $B$, is the quasi-affine variety $K_{\C}/N$.
This stratification agrees with the symplectic stratification of section
\ref{symplectic}.
In particular, using the Iwasawa decomposition $K_{\C} = KAN$, we may view the
top stratum as $KA$, the open subset of the implosion
corresponding to the interior of the positive Weyl chamber
for $K$.

The simplest example
as discussed in \cite{GJS:implosion}, is $K= SU(2)$. Now the $N$ action on $K_{\C}=SL(2, \C)$ is:
\[
\left(\begin{array}{cc}
x_{11} & x_{12} \\
x_{21} & x_{22}
\end{array} \right) \mapsto
\left(\begin{array}{cc}
x_{11} & x_{12} \\
x_{21} & x_{22}
\end{array} \right) 
\left(\begin{array}{cc}
1 & n \\
0 & 1
\end{array} \right) =
\left(\begin{array}{cc}
x_{11} & x_{12} + n x_{11}\\
x_{21} & x_{22} + n x_{21}
\end{array} \right)
\]
with invariant ring freely generated by $x_{11}$ and $x_{21}$,
so $K_{\C}/\!/N = \C^2$. There are two strata, the top one
$SL(2,\C)/N = \C^2 - \{ 0 \}$ and the bottom one $\{ 0 \}$.
(As the closed Weyl chamber for $SU(2)$ is $[0, \infty)$, these
coincide with the symplectic strata $SU(2) \times (0, \infty)$
and $(SU(2) \times \{0 \})/SU(2)$).

So we see, as in the general case, that the implosion provides
an affine completion of the quasi-affine top stratum. Notice that the
canonical morphism $K_{\C} \rightarrow K_{\C}/\!/N = \C^2$ defined by
$\left( \begin{array}{cc}
x_{11} & x_{12} \\
x_{21} & x_{22}
\end{array} \right) \mapsto 
\left( \begin{array}{c}
x_{11} \\
x_{21}
\end{array} \right)$ is not surjective, but instead 
has image the constructible set $\C^2 -\{ 0 \}$.

In this case the strata actually fit together to form a smooth variety,
but if $K$ has a simple factor of rank greater than one, the implosion
is always singular.

\medskip
This picture has been generalised by Kirwan \cite{K} to the case of
quotients $K_{\C}/\!/U_P$ where $U_P$ is the unipotent radical of
a parabolic subgroup $P$. This nonreductive quotient still exists
as a variety and there is an interpretation in terms of a generalised
version of the symplectic implosion construction of section
\ref{symplectic}. These spaces are referred to as partial symplectic
implosions. They have an action of $K_{\C} \times L_P$ where $L_P$
is the reductive Levi subgroup of $P$ (recall $P$ is the semidirect product
 $U_P  \rtimes L_P$). 

\section{Hyperk\"ahler implosion}
\label{sec:hk-implosion}

In \cite{DKS}  we considered an analogue of the universal implosion for
hyperk\"ahler geometry. The starting point is the observation by 
Kronheimer \cite{Kronheimer:cotangent} that $T^*K_{\C}$ carries a complete hyperk\"ahler metric
that is preserved by an action of $K \times K$. This action is
not only isometric but also triholomorphic, that is, it
is preserves each individual complex structure $\sf I, \sf J,\sf K$.

Kronheimer's construction proceeds by identifying $T^* K_{\C}$ with
the moduli space of solutions to Nahm's equations
\[
\frac{dT_i}{dt } + [T_0, T_i] = [T_j, T_k],
\]
where $(ijk)$ is a cyclic permutation of $(123)$, for smooth maps $T_i : [0,1] \rightarrow
\kf$. The moduli space is formed by quotienting by the gauge group
of maps $g : [0,1] \rightarrow K$ such that $g(0)= g(1) = Id$. 

The residual gauge action by gauge transformations not necessarily
equal to the identity at the endpoints
$0,1$ gives rise to the hyperk\"ahler $K \times K$ action. 

Note also that there is an isometric $SO(3)$ action given by
rotating the triple $(T_1, T_2, T_3)$ of Nahm matrices. This action
is {\em not} triholomorphic but acts transitively on the 2-sphere
of complex structures.

The identification of the Kronheimer moduli space with $T^*K_{\C}$
involves of course a choice of complex structure $\sf I$. However all
such complex structures are equivalent under the $SO(3)$ action.
Note also that the $\sf I$-holomorphic symplectic structure defined by the
holomorphic parallel 2-form $\omega_{\sf J} + i \omega_{\sf K}$ is just the
standard $K_{\C} \times K_{\C}$-invariant holomorphic symplectic structure that $T^*K_{\C}$ has as the cotangent bundle of a complex manifold. (We shall usually use the term{\em complex-symplectic  structure} for
holomorphic symplectic structure in this paper).

 $T^*K_{\C}$ is thus the hyperk\"ahler analogue of the symplectic
$K \times K$-space $T^*K = K_{\C}$.
As the universal symplectic implosion is the nonreductive quotient
$K_{\C} /\!/N$,
it makes sense in the hyperk\"ahler setting to consider a suitable
reduction of $T^*K_{\C}$ by $N$, more specifically the
complex-symplectic quotient (in the sense of geometric invariant theory) of $T^*K_{\C}$ by $N$.

As the complex-symplectic structure on $T^*K_{\C}$ is the standard one,
its associated moment map is just projection onto the $\kf^*_{\C}$ factor
of $T^*K_{\C} = K_{\C} \times \kf^*_{\C}$. The zero locus for this
moment map is therefore $K_\C \times \n^\circ$ where
\( \n^\circ \) is the annihilator in \( \kf_{\C}^* \) of the Lie algebra
\( \n \) of $N$. 

We are therefore led to define the universal hyperk\"ahler implosion
for $K$ to be the geometric invariant theory (GIT) quotient
\( (K_\C \times \n^\circ) \symp N \) where \( N \) is a 
maximal unipotent subgroup of the complexified group \( K_\C \).
It is sometimes convenient to choose an invariant
inner product, and identify the annihilator $\n^\circ$ with
the opposite Borel subalgebra $\bmf$).

As $N$ is nonreductive, it is a nontrivial result that  the algebra of
$N$-invariants is finitely generated and hence 
the quotient exists as an affine variety.
This was shown in the case $K=SU(n)$ in \cite{DKS} and in general
follows from results of Ginzburg-Riche \cite{GinzburgRiche}
(see the discussion in \cite{DKS-Arb} for example). 

The universal hyperk\"ahler implosion 
carries a complex-symplectic action of \( K_\C \times T_\C \) where  \( T \)
is the standard maximal torus of $K$. The $K_{\C}$ action is just
left translation on the $K_{\C}$ factor, while the the $T_{\C}$ action is
right translation on the $K_{\C}$ factor together with the adjoint action
on the $\n^\circ$ factor. Of course the fact we are restricting to $\n^\circ$
means that the right $K_{\C}$ action on $K_{\C} \times \kf_{\C}^*$ is broken
to a $T_{\C}$ action.

A naive guess might be that, by analogy with the symplectic case,
the complex-symplectic torus reductions of the implosion will
give us the coadjoint orbits for the complex Lie algebra $\kf_{\C}$.
However this cannot be exactly right, as only semisimple coadjoint orbits
in the complex Lie algebra are closed.
The complex-symplectic quotients
by the torus action are instead
the Kostant varieties; that is, the varieties in \( \kf_{\C}^* \)
obtained by fixing the values of the invariant polynomials for this
Lie algebra \cite{Chriss-G:representation, Kostant:polynomial}.  
The Kostant varieties are in general stratified
spaces whose strata are distinct complex coadjoint orbits. The
minimal stratum is the semisimple orbit and the top stratum
 is the regular orbit,
 which is open and dense in the Kostant variety with complement
of codimension at least 2. (For $K_{\C} = SL(n, \C)$ the elements
of the regular orbit are characterised by the minimal polynomial
being equal to the characteristic polynomial, the latter being
fixed by the choice of Kostant variety).

Note that, just as the symplectic implosion has real dimension
$\dim_{\R}K+ {\rm rank \;}K$, so the hyperk\"ahler implosion
has {\em complex} dimension equal to $\dim_{\C} K_{\C} + {\rm rank} \; K_{\C}$,
consistent with the fact that the Kostant varieties have complex dimension
$\dim_{\C} K_{\C} - {\rm rank \;} K_{\C}$.

%One may, by analogy with the symplectic case, then
%implode a general space with hyperk\"ahler $SU(n)$ action by taking its
%product with the universal implosion and performing the hyperk\"ahler
%reduction by the diagonal $K$ action.
%For general compact groups a direct construction
%of a hyperk\"ahler metric on the nonreductive quotient 
%\( (K_\C \times \n^0) \symp N \) is not yet available, although in \cite{DKR} we
%gave an alternative approach to hyperk\"ahler implosion

\section{Hyperk\"ahler quiver diagrams}
\label{sec:hk-quiver}

The description in the previous section is rather abstract and although
it makes plain the complex-symplectic structure, it is less clear
that this actually comes from a hyperk\"ahler metric.

In \cite{DKS} we considered the case when \( K=\SU(n) \).
In this situation the universal hyperk\"ahler implosion can be identified with 
a hyperk\"ahler quotient using quiver diagrams, and thus
can be seen to be genuinely a stratified hyperk\"ahler space rather than just a 
complex-symplectic one. 

We shall consider quivers $Q= (Q_0, Q_1)$ where $Q_0$ is the set of vertices
and $Q_1$ the set of edges. For each edge $e \in Q_1$, we denote 
$o(e)$ and $i(e)$ the outgoing and incoming vertices of the edge.
To each vertex $j$ we associate a complex vector space $V_j$ of dimension
$N_j$.

In the simplest case one can associate to the quiver the flat
quaternionic space
\[
M = \oplus_{e \in Q_1} \hom (V_{i(e)}, V_{o(e)}) \oplus
\hom (V_{o(e)}, V_{i(e)})
\]
and the group $K = \prod_{j \in Q_0} U(V_j)$, with its natural action on $M$:
\begin{equation*}
  \alpha_e \mapsto g_{o(e)} \alpha_e g_{i(e)}^{-1},\quad
  \beta_e \mapsto g_{i(e)} \beta_e g_{o(e)}^{-1} \qquad (e \in Q_1),
\end{equation*}
In more physical language, to each edge joining vertices
labelled by dimensions $N_i$ and $N_j$ we associate the hypermultiplets
$\HH^{N_i N_j}$ transforming in the bifundamental representation of $U(N_i) \times U(N_j)$. Fixing a complex structure and identifying this with $\hom (\C^{N_i}, \C^{N_j}) \oplus \hom(\C^{N_j}, \C^{N_i})$ as above corresponds physically
to decomposing the hypermultiplet into chiral and antichiral multiplets.

This action preserves the hyperk\"ahler structure so one may
form the hyperk\"ahler reduction $M \hkq K$. More generally, 
one may hyperk\"ahler reduce by a subgroup $K_1$ of $K$, so that the
quotient $M \hkq K_1$ retains a residual hyperk\"ahler action
of ${\sf N}_{K}(K_1)/K_1$ where ${\sf N}_{K}(K_1)$ denotes the normaliser of
$K_1$ in $K$. In particular, one may define a normal subgroup
$K_1$ of $K$ by choosing a subset $Q \subset Q_0$ and defining
$K_1 = K_Q := \prod_{j \in Q} U(V_j)$. That is, we `turn off' the
action at the nodes in $Q_0 - Q$. The hyperk\"ahler quotient
now has a residual action of $K/K_Q \cong \prod_{j \notin Q} U(V_j)$.

The vertices $j \in Q$ where the group still acts are called 
{\em gauge nodes} and the vertices $j \in Q_0 - Q$ where the
action has been turned off are the {\em flavour nodes}. The
gauge nodes are denoted by circles and the flavour nodes by square boxes.

\begin{example} \label{nilpotentquiver}

Consider the $A_{n}$ diagram with dimension vector $(1,2, \ldots,n)$
where the $n$-dimensional node is a flavour node. (The figure shows the $n=6$ case)

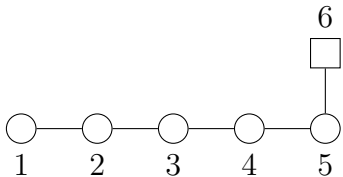
\begin{figure}
	[t] 
\centering 
	\begin{tikzpicture}
%		\node (dg1) at (-1,0) [gauge, label=below:{1}]{}; 
%		\node (dg2) at (0,0) [gauge, label=below:{2}]{}; 
%		\node (dg3) at (.7,.7) [gauge, label=below:{1}]{};
%		\node (dg4) at (.7,-.7) [gauge, label=below:{1}]{}; 
%		\node (df2) at (0,1) [flavor, label=above:{1}]{}; 
%		\draw (dg1)--(dg2); 
%		\draw (dg2)--(dg3); 
%		\draw (dg2)--(dg4); 
%		\draw (dg2)--(df2); 
%		
%		\node (bg1) at (5.5,0) [gauge, label=below:{\symmlabel{1}}]{}; 
%		\node (bg2) at (6.5,0) [gauge, label=below:{\symmlabel{2}}]{}; 
%		\node (bg3) at (7.5,0) [gauge, label=below:{\symmlabel{[1]\wr S_2}}]{}; 
%		\node (bf2) at (6.5,1) [flavor, label=above:{1}]{}; 
%		\draw (bg1)--(bg2)--(bg3); 
%		\draw (bg2)--(bf2); 
		
%		\draw[->,thick] (3,0) -- (4,0);
		
		\node (ag1) at (-2,-3.5) [gauge, label=below:{1}]{}; 
		\node (ag2) at (-1,-3.5) [gauge, label=below:{2}]{}; 
		\node (ag3) at (0,-3.5) [gauge, label=below:{3}]{}; 
		\node (ag4) at (1,-3.5) [gauge, label=below:{4}]{}; 
		\node (ag5) at (2,-3.5)  [gauge, label=below:{5}]{}; 
%		\node (af1) at (-2,-2.5) [flavor, label=above:{0}]{}; 
		\node (af5) at (2,-2.5) [flavor, label=above:{6}]{}; 
		\draw (ag1)--(ag2)--(ag3)--(ag4)--(ag5); 
%		\draw (ag1)--(af1); 
		\draw (ag5)--(af5); 
		
%		\node (cg1) at (6,-3.5) [gauge, label=below:{\symmlabel{1}}]{}; 
%		\node (cg2) at (7,-3.5) [gauge, label=below:{\symmlabel{\qquad 1]\wr S_2}}]{}; 
%		\node (cg3) at (8.3,-3.5) [gauge, label=below:{\symmlabel{1}}]{}; 
%		\node (cf1) at (5,-3.5) [flavor, label=below:{\symmlabel{[1}}]{}; 
%		\draw (cf1)--(cg1)--(cg2)--(cg3); 
		
%		\draw[->, thick] (3,-3.5) -- (4,-3.5);
	\end{tikzpicture}
	\caption{Quiver for the nilpotent cone of $A_5$.}
\end{figure}

So we hyperk\"ahler reduce by $U(1) \times \ldots \times U(n-1)$
and leave the residual action $U(n)$. The hyperk\"ahler quotient
is known by the work of Kobak-Swann \cite{Kobak-S:finite} 
(see also \cite{KP}) to be the nilpotent
variety for $A_{n-1}=SL(n, \C)$

\end{example}

This motivated the quiver description of hyperk\"ahler implosion
for $K = SU(n)$ developed in \cite{DKS}. The implosion is required
to have a $SU(n) \times T$ action with hyperk\"ahler reduction by $T$
giving the Kostant varieties, in particular reduction at level zero
giving the nilpotent variety. It is natural therefore to consider
the same quiver as above, but with the action of
 \( H = \prod_{j=1}^{n-1}\SU(j) \), rather than $K = \prod_{j=1}^{n-1}
U(j)$. The resulting hyperk\"ahler quotient
$M \hkq H$  is a stratified hyperk\"ahler space
 with a residual action of
the torus \( T = K / H \) 
as well as a commuting action of \( \SU(n) \).

We can also consider the implosion as a complex-symplectic
quotient. It is the geometric invariant theory quotient, 
of the zero locus of the complex moment map \( \mu_{\C} \) for the
\( H \) action, by the 
complexification
\begin{equation*}
  H_\C  = \prod_{j=1}^{n-1}\SL(j,\C)
\end{equation*}
of \( H \), 

The complex moment map equation \( \mu_\C =0 \) is equivalent to the
equations
\begin{equation}
  \label{eq:mmcomplex}
   \beta_{i+1} \alpha_{i+1} - \alpha_{i}\beta_{i}  = \lambda^\C_{i+1} I \qquad
  (i=0,\dots, n-2),
\end{equation}
for (free) complex scalars \( \lambda^\C_1,\dots,\lambda^\C_{n-1} \).
The complex numbers \( \lambda_i \) combine to give the
complex-symplectic moment map for the
 residual action of \( K_\C/H_\C \) which we can
identify with the maximal torus \( T_\C \) of
\( K_\C \).

%It is often useful to consider the endomorphism
%\begin{equation*}
%  X = \alpha_{r-1} \beta_{r-1} \in \Hom (\C^n,\C^n),
%\end{equation*}
%which is invariant under the action of \( \tH_\C \) and transforms by
%conjugation under the residual \( \SL(n,\C) \) action.  

Note that, as usual with linear hyperk\"ahler quotients at level zero, we also
have an $Sp(1)$ action on the implosion that {\em rotates} the
complex structures. If we view the quaternionic summands 
$\hom (V_{i}, V_{i+1}) \oplus \hom (V_{i+1}, V_{i})$ associated to each edge of the quiver as quaternionic space $\HH^{N_i N_{i+1}}$  then the quiver group $H$ may
be viewed as acting on $\HH^{N_i N_{i+1}}$ on the left while the
quaternionic structure is acting on the right by $-i,-j,-k$ etc.
Now multiplication by unit quaternions on the right gives an isometric
action, rotating complex structures, and commuting with the action of $H$.
It therefore acts on the hyperk\"ahler moment map
$\mu : M \rightarrow \h^* \otimes \R^3$ by rotation on $\R^3$ and hence
preserves the hyperk\"ahler quotient at level zero. Moreover, as the
level is zero and the moment map is homogeneous quadratic, we
have a scaling action of the positive reals. We can summarise this
as saying the SU(n)-implosion has a conical structure, and as such
is expected to fit into the symplectic duality framework discussed
in section \ref{sec:duality}.

\medskip
For other classical groups we do not as yet have a quiver description
of the implosion. This is because the analogues of the quiver description
of the nilpotent varieties involve {\em orthosymplectic quivers}, that is,
quivers where the groups attached to the vertices are alternately orthogonal
and symplectic groups \cite{Kobak-S:finite}. Unlike the unitary groups, we cannot write these groups
as extensions of tori by subgroups, so we cannot mimic the above construction
by considering quivers with just the subgroups acting.

\section{Moore-Tachikawa category} \label{sec:MT}

In \cite{MT} Moore and Tachikawa proposed a category 
whose objects were complex semisimple or reductive groups and where morphisms
between $G_1$ and $G_2$ are complex-symplectic manifolds with
$G_1 \times G_2$ action. (Strictly speaking a morphism is a triple
$(X, G_1, G_2)$ where $X$ is such a complex-symplectic
manifold, ie the ordering of the objects is specified). There is also supposed to be a 
commuting circle action acting on the complex-symplectic form with
weight 2. Composition of morphisms $X \in {\rm Mor}(G_1, G_2)$
and $Y \in {\rm Mor}(G_2, G_3)$ proceeds by forming the product
$X \times Y$ with $G_1 \times (G_2 \times G_2) \times G_3$ action
and then taking the complex-symplectic quotient by the diagonal
$G_2$ action. The resulting quotient is complex-symplectic with
residual $G_1 \times G_3$ action so lies in Mor$(G_1, G_3)$
as required. The Kronheimer space $T^*K_{\C}$ is complex-symplectic
with $K_{\C} \times K_{\C}$ action and defines the identity element
in Mor$(G,G)$ with $G = K_\C$.

In this picture the implosion for $K$ may be viewed as an
element of Mor$(K_{\C}, T_{\C})$. The process of imploding
a complex-symplectic manifold with $K_{\C}$ action to obtain
a manifold with $T_{\C}$ action, as described in
section \ref{symplectic} but in the complex-symplectic case,
is now exactly that of composition of morphisms with the implosion, to
obtain a map:
\[
{\rm Mor}(1, K_{\C}) \rightarrow {\rm Mor}(1, T_\C)
\]
Note that one could enrich the data of complex-symplectic manifolds to
 hyperk\"ahler manifolds in these
definitions, using the fact that the complex-symplectic quotient
by $G_2$ coincides with the hyperk\"ahler quotient by the maximal compact
subgroup of $G_2$. 
However now $T^*K_{\C}$ is no longer exactly the identity, as pointed out by Moore-Tachikawa. The metric is shifted by a factor representing the
length of the interval on which the Nahm data is defined to produce the
Kronheimer space.

\section{Symplectic duality}
\label{sec:duality}

It is conjectured that there is a 
duality between certain complex-symplectic (that is, holomorphic symplectic) varieties, that
physically may be interpreted as duality (the notion of duality is explained below and is different than other forms of dualities in physics) between Higgs and Coulomb branches of
a 3d $N=4$ theory

The complex-symplectic varieties concerned usually in fact have a
hyperk\"ahler structure, and arise either as hyperk\"ahler cones or as
deformations thereof. In many cases the Higgs branch cone occurs as
the zero level set of a hyperk\"ahler quotient construction $M \hkq G$
(the moduli space of vacua), while the deformations occur by moving
the level set away from zero. In physics the resulting deformation
parameters are called Fayet-Iliopoulos parameters.

For symplectic duality
constructions we want the complex-symplectic varieties to have a
circle action that acts on the complex symplectic form with weight 2
(in terms of the hyperk\"ahler structure, the circle action fixes one
complex structure $\sf I$ but rotates the $\sf J,{\sf K}$ so the $\sf I$-holomorphic
form $\omega_{\sf J} + i \omega_{\sf K}$ is scaled rather than being invariant
under the action).

As mentioned in \S \ref{sec:hk-quiver}, linear hyperk\"ahler quotients
at level zero have a $Sp(1)$ action  rotating the complex structures.
Making a deformation that breaks this $Sp(1)$ down to the circle action
fixing the specific complex structure $\sf I$ corresponds to changing the
level set to $(\lambda,0,0)$ where $\lambda \in \mathfrak g^*$
As the level set at which the hyperk\"ahler reduction is performed must lie in the centre of $G$, the number of deformation parameters 
is the dimension of the center of $G$.

On the Coulomb side, the deformation parameters are the masses.
The duality is supposed to interchange the rank of the hyperk\"ahler isometry
group of a space and the number of deformation parameters for its dual.
More precisely, the Cartan algebra of the flavour group of the Higgs branch
is identified with the space of mass parameters, and the Cartan algebra of the
flavour group of the  Coulomb branch with the space of
Fayet-Iliopoulos parameters.

Nakajima (see \cite{N1} for example) has suggested that in the case when
the Higgs branch is a hyperk\"ahler quotient $M \hkq G$ by a compact group
$G$, the Coulomb branch should be birational to
$T^* (T_{\C}^{\vee})/W$, the quotient by the Weyl group of the cotangent bundle of the complexified dual maximal torus
of $G$. We therefore expect
\[
  \dim_{\R} (\rm {Coulomb \; branch}) = 4 \; {\rm rank \;} G.
 \]
 Physically, the birational equivalence represents quantum corrections to the
 classical description of the Coulomb branch.
 
\bigskip
One example where the theory is completely worked out is {\em hypertoric manifolds}, that is, hyperk\"ahler quotients of flat quaternionic space by
tori. (See \cite{BLP} for example).
As in \cite{BD} we consider quotients of $\HH^d$ by a subtorus $N $ of $T^d$.
The torus is defined by vectors $u_1, \ldots, u_d \in \R^n$ : explicitly
we define $\n = {\rm Lie} N$ to be the kernel of the map 
$\beta : \R^d = {\rm Lie}
T^d \rightarrow \R^n$ defined by $\beta : e_i \mapsto u_i$, where
$e_1, \ldots, e_d$ is the standard basis for $\R^d$. On the Lie
algebra level, we have an exact sequence
\[
0 \rightarrow \n \rightarrow \R^d \overset \beta \rightarrow \R^n \rightarrow 0
\]
On the Lie group level we have:
\[
1 \rightarrow N \rightarrow T^d \rightarrow T^n \rightarrow 1
\]
The  hypertoric $M = \HH^d \hkq N$ has real dimension $4d - 4(d-n) = 4n$.
and admits a residual action of the quotient torus $T^n = T^d / N$.
The number of deformation parameters for $M$ is rank $N = d-n$
and the rank of the isometry group is $n$. 

Now the dual hypertoric variety is defined to be the hyperk\"ahler quotient
of $\HH^d$ by the dual torus $\hat{T^n}$
\[
1 \rightarrow \hat{T}^n \rightarrow \hat{T}^d \rightarrow \hat{N} \rightarrow 1
\]
Now the number of deformation parameters is $n$ and the rank of the isometry group is ${\rm rank} \; \hat{N} = {\rank N} = d-n$, in accordance with the
principle of symplectic duality. The dimension of the dual hypertoric
is $4(d-n)$, illustrating how dimension can change under duality.

As usual in toric or hypertoric geometry, this duality can be viewed
as a combinatorial phenomenon, in this case known as Gale duality.
Given a vector space $V$ of dimension $n$ with spanning 
vectors $u_1, \ldots, u_d$, we can form the space of linear dependency relations
$\{ (\alpha_1, \ldots, \alpha_d) : \sum_{i=1}^{d} \alpha_i u_i =0 \}$.
This is a $d-n$ dimensional vector space $W$ with $d$ distinguished elements
$w_1, \ldots, w_d$ in the dual vector space $W^*$ defined by 
$w_i : (\alpha_i, \ldots, \alpha_d) \mapsto \alpha_i$. This 
duality, interchangeing $n$ and $d-n$, implements
the above duality between the hypertorics of dimension $4n$ and $4(d-n)$.

In this case, both the Higgs and Coulomb branches are given by
finite-dimensional hyperk\"ahler quotients. However there are cases
where one space is given by such a construction but its dual is not--we call these {\em non-Lagrangian} theories.

\bigskip
Various relations between a quiver variety and its symplectic dual
have been developed in the physics literature.
 
The crucial concept here is that of a {\em balanced} node. In the
case of a unitary quiver with dimensions $N_j$ at nodes $j$, the
{\em balance} of a node $j$ is
\[
-2 N_j + \sum_{k \; {\rm adjacent \; to} \; j}  N_k
\]
and we say the node is {\em balanced} if the balance is zero.

For a nice physical theory we would like all the gauge nodes to have
balance greater than or equal to $-1$. If this
holds and there is a node with balance equal to $-1$ the quiver is
called {\em minimally unbalanced}, while if all nodes have nonnegative
balance with at least one of positive balance, we say it is {\em positively
  balanced}.

In the case of unitay quivers, the balanced gauge nodes should form the Dynkin diagram of the
semisimple part of (a subgroup of) the hyperk\"ahler isometry group of the dual space.
(Unbalanced nodes give abelian symmetries). This refines the earlier
idea that deformation parameters coming from the unitary gauge nodes
should give an abelian algebra of symmetries in the dual--if the nodes are balanced
then the associated abelian symmetry group is realised as the maximal
torus of a larger semisimple group.

For example, in the nilpotent variety quiver of Example \ref{nilpotentquiver} 
all nodes except the final flavour node are balanced. This gives
an $A_{n-1}$ Dynkin diagram which should give $SU(n)$ symmetry group of
the dual. In fact the dual is still the nilpotent variety.

\begin{example}
Consider the quiver diagram in Figure 2 corresponding to the hyperk\"ahler quotient
$\HH^d \hkq U(1)$

%(diagram)
\begin{figure}
	[t] 
\centering 
	\begin{tikzpicture}		
		\node (ag5) at (2,-3.5)  [gauge, label=below:{1}]{}; 
		\node (af5) at (2,-2.5) [flavor, label=above:{$d$}]{}; 
		\draw (ag5)--(af5); 
		
]	\end{tikzpicture}
	\caption{$U(1)$ with $d$ flavors.}
\end{figure}
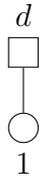

where we have 1 gauge node (with dimension 1) and 1 flavour node 
(with dimension $d$)

This is a hypertoric, with symplectic dual $\HH^d \hkq T^{d-1}$. The latter
space gives the cyclic Kleinian singularity $\C^2/\Z_d$ or its deformations,
the $A_{d-1}$ multi-instanton metrics whose topology is generated by a chain
of $d-1$ rational curves with self-intersection $-2$.

If $d \neq 2$ than we have no balanced nodes in the diagram, but if
$d=2$ then the gauge node is balanced. This reflects the fact that
for $d=2$ the dual space is Eguchi-Hanson which has a triholomorphic $SU(2)$
action, an enlargement of the triholomorphic $U(1)$ action that occurs
for general $d$.
\end{example}

One can study the varieties occuring in symplectic duality by
finding the Hilbert series of their coordinate ring (the chiral ring
in physics terminology). This series counts the dimension $m_d$ of the degree
$d$ parts of the ring
\[
HS(t) =  \sum_{d=0}^{\infty} m_d t^d.
\]
The variable $t$ is called the {\em fugacity}.

Cremonesi-Hanany-Zaffironi \cite{CHZ} have derived a formula, the {\em monopole formula} to compute the Hilbert series of the Coulomb branch of a quiver
variety obtained as a hyperk\"ahler reduction of a flat quaternionic space
by a group $G$. We are counting monopole operators whose gauge field has
a Dirac monopole singularity, with associated magnetic charge living in the
weight lattice $\Gamma_{\hat{G}}$ of the Langlands dual $\hat{G}$.
Their formula involves contributions from the stabiliser groups
of each element of the lattice:
%\[
\begin{equation}
HS(t) = \sum_{m \in \Gamma_{\hat{G}}/W_{\hat{G}}} t^{2\Delta(m)} P_{G}(m,t)
\label{HS1}
\end{equation}
%\]
Here
\[
  P_G(m,t) = \prod_{i} \frac{1}{1- t^{2d_i(m)}}
\]
where the $d_i(m)$ are the exponents of the stabiliser
group $G_m = {\rm Stab}_{G}(m)$--that is, the degrees of the generators
(Casimirs) for the ring of invariants of $G_m$ under the adjoint representation.
We can also interpret $P_G(m,t)$ as the Poincar\'e polynomial of the 
classifying space $BG_m$.

The term $\Delta(m)$ is given by
\[
\Delta(m) = - \sum_{\alpha \in R^+} | \alpha(m) | + \frac{1}{2} \sum_{b} |b(m)|
\]
where $R^+$ denotes the set of positive roots in $G$, and the second sum
is taken over the weights in the given representation $M$.

Plethystic techniques have been developed (eg \cite{FHH}) to
compute from the Hilbert series the generators, relations and higher-order
syzygies of the chiral ring. 

Note that the $t^2$ term of the Hilbert series is expected to give the
dimension of the global symmetry group.

\section{Duals of implosions}
\label{sec: dual-implosion}

We now consider what kind of space would be dual to the $SU(n)$ implosion.
The latter space has an action of $SU(n) \times T^{n-1}$, so this suggests
we look at a quiver whose balanced nodes give the Dynkin diagram $A_{n-1}$
and whose unbalanced nodes give the torus factor.

We consider the quiver from Example \ref{nilpotentquiver} that gives the nilpotent variety.
Now replace the flavour node (box) with dimension $n$ by a bouquet of $n$ $U(1)$ nodes attached to the $(n-1)$-dimensional gauge node. This ensures that the
$(n-1)$-dimensional node remains balanced, as well as the gauge nodes
lower down the chain. So the balanced nodes do form the $A_{n-1}$ Dynkin diagram
as required, giving a $SU(n)$ action on the dual. The $U(1)$ nodes are unbalanced (for $n \neq 3$) and generate a $T^{n-1}$ action on the dual (it is $T^{n-1}$,
 not $T^n$, as one $U(1)$ `decouples', ie 
acts trivially. This is the diagonally embedded 
$U(1) \hookrightarrow U(1)^n \times \prod_{k=1}^{n-1} Z(U(k))$, where $Z$ denotes
the centre). Note that the balance of the $U(1)$ nodes is always at least
$-1$, and is positive for $n \geq 4$.

\begin{example}
If $n=2$ this is just an $A_3$ diagram with dimension 1 at each
node. As the diagonal $U(1)$ acts trivially this represents the trivial
hypertoric
$\HH^2 \hkq T^2$ and its dual is $\HH^2 \hkq \{ 1 \} = \HH^2$. This is
correct as the universal hyperk\"ahler implosion for $SU(2)$ is indeed
$\HH^2$.
\end{example}

\begin{example}
  If $n=3$ we have a star-shaped quiver (affine $\tilde{D}_4$ Dynkin diagram)
  with dimension 2 at the central node
  and dimension 1 at the four nodes radially connected to it (one from the
  tail of the truncated $A_2$ diagram and three from the bouquet). Uniquely
  in this case all nodes (even the bouquet ones) are balanced, so we expect,
  after decoupling, an $SO(8)$ symmetry in its dual.

  This is correct, as the $SU(3)$ universal hyperk\"ahler implosion
  may be identified with the Swann bundle of the quaternionic K\"ahler
  Grassmannian $\tilde{\rm Gr}_4(\R^8) = SO(8)/S(O(4) \times O(4))$
  of oriented 4-planes in $\R^8$. The $SO(8)$ symmetry of the quaternionic
  K\"ahler base lifts to a symmetry of the hyperk\"ahler Swann bundle
  (see Example 8.7 of \cite{DKS} for a discussion).
  \end{example}
  
As the $SU(n)$ implosion has been described as a reduction by a product
of special unitary groups in \ref{sec:hk-quiver}, we expect it has no
deformation parameters. This checks with the fact that the proposed dual
has no residual hyperk\"ahler isometries, as all nodes are gauge and not flavour
nodes. 

In fact, we expect for general groups that the implosion has no
 deformation parameters, as we obtain it as the nonreductive quotient
 $(K_{\C} \times \n^\circ ) \symp N$ and the maximal unipotent group $N$ has
trivial maximal torus so no characters.

For a global symmetry of $SU(n) \times U(1)^{n-1}$ we expect the coefficient of the $t^2$ term in the Hilbert series to be $n^2-1+n-1=n^2+n-2.$ In addition, due to the balance of $n-3$ of each U(1) node in the Bouquet, there are generators of the chiral ring which arise from the U(1) nodes that contribute 2 per each U(1) at order $t^{n-1}$. These correspond to one monopole operator of positive charge and one with negative charge under the corresponding U(1) global symmetry. We expect the Hilbert series to get contributions
$$HS_n = (n^2+n-2)t^2 + 2nt^{n-1}+\ldots$$
Let us see how this fits in examples.
For $n=2$ we get $4t$ that represent the 4 generators of $\HH^2$. They contribute 6 more quadratic bilinears that enhance the global symmetry from $SU(2) \times U(1)$ to $Sp(2)$.
For $n=3$ the affine $D_4$ quiver indeed confirms that the global symmetry is enhanced from $SU(3) \times U(1)^2$ to $SO(8)$. For $n>3$ perturbative computations confirm the $t^2$ coefficient.

One can further refine the expression for the Hilbert Series in equation (\ref{HS1}) by introducing a fugacity $z_i$ for each magnetic charge $m_i$ of $U(1)_i$ in the bouquet for $i=1\ldots n$, resulting in a function of $n+1$ variables $HS(t,z_i)$. This expression can be further integrated
\[
(1-t^2)^{n-1}\prod_i \oint_{|z_i|=1} \frac{dz_i}{z_i}HS(t,z_i)
\]
resulting in the expression for the Hilbert series of the nilpotent cone of $SL(n)$ which takes a particularly simple form
\[
\frac{\prod_{i=1}^n (1-t^{2i})}{(1-t^2)^{n^2}}
\]
This constitutes a non trivial test of the proposed quiver for the $SU(n)$ implosion.

\bigskip
We can also check that this is consistent with Nakajima's picture.

The rank of the group $U(1)^{n-1} \times \prod_{i=1}^{n-1} U(i)$
by which we quotient in the bouquet quiver is
$\frac{1}{2}(n+2)(n-1)= \frac{1}{2}(n^2+n-2)$ and the real dimension of the implosion is
\[
  \dim_{\R} SL(n,\C) + \dim_{\R} (T_{\C}^{n-1}) = 2(n^2 + n-2).
\]
Going in the reverse direction, the implosion is produced as a hyperk\"ahler
quotient by $\prod_{i=1}^{n-1} SU(i)$ which has rank $\frac{1}{2}(n-1)(n-2)$.
The quaternionic dimension of the bouquet quiver variety is
\[
 n(n-1) + \sum_{i=1}^{n-2} i(i+1)  - (n-1 + \sum_{i=1}^{n-1} i^2)
\]
which works out as $\frac{1}{2}(n-1)(n-2)$ as desired.

For example,
if $n=3$ then we have the affine $\tilde{D}_4$ Dynkin diagram,
giving one of Kronheimer's examples  \cite{Kronheimer:ALE} of real dimension 4, ie quaternionic
dimension 1. This corresponds to the fact that the $SU(3)$ implosion
is a hyperk\"ahler quotient of a linear space by $SU(2)$.

\medskip
We also make some remarks on partial hyperk\"ahler implosions,
ie complex symplectic quotients of $T^*K_{\C}$ by the unipotent radical
$U_P$ of a parabolic $P$. (It as as yet a conjecture that these exist as
algebraic varieties, that is, that the algebra of $U_P$-invariants
in $K_\C \times \mathfrak{u}_{P}^{\circ}$ is finitely generated).

In the case $K= SU(n)$, of course, the
parabolics are indexed by ordered partitions $n =n_1 +  \ldots+ n_r$
and the corresponding Levi subgroup is $S( GL(n_1, \C) \times \ldots \times
GL(n_r,\C))$.

As $SL(n,\C)/P = SU(n)/S(U(n_1) \times \ldots \times U(n_r))$, we see
that
\[
  \dim_{\R} P = n^2 -2 + \sum_{i=1}^{r} n_i^2
\]
and
\[
\dim_{\R} U_P = n^2 - \sum_{i=1}^{r} n_i^2
\]
so the dimension of the partial implosion should be
\[
  \dim_{\R} (SL(n,\C) \times \mathfrak{u}_{P}^{\circ}) /\!/ U_P
  = 2 (n^2 -2 + \sum_{i=1}^{r} n_i^2)
\]
Note that as $\sum_{i=1}^{r} n_i =n$, the
sum $\sum_{i=1}^{r} n_i^2$ has the same parity as $n^2$ so the expression
inside the bracket above is even, as required.

If all $n_i=1$ of course $P$ is the Borel and we recover the dimension of
the standard implosion as above.

\medskip
A natural candidate for the dual would be the quiver diagram we obtain
by taking the basic diagram for the nilpotent quiver, excising the dimension
$n$ flavour node, and then attaching $r$ legs, each of them
an $A_{n_i}$ quiver with the dimension $n_i$ node next to the dimension
$n-1$ node of the original diagram.

So the remaining nodes of the original diagram are all balanced, giving
an $SU(n)$ symmetry in the implosion.
Moreover on each leg, all nodes except the $\dim n_i$ ones are balanced,
yielding $SU(n_i)$ symmetries for $i=1,...,r$.
Also, the $r$ unbalanced nodes (ie the $\dim n_i$ ones of the attached legs)
would yield, after decoupling,  $r-1$ Abelian symmetries. These nodes
have balance $n - n_i -2$ which is always at least $-1$ and is positive
unless our partition is $n = (n-2) + 2, (n-2) + 1+ 1$ or $(n-1) +1$.

So overall, we would get $SU(n) \times S(U(n_1) \times...\times U(n_r))$ symmetry, as required.

The group by which we perform the hyperk\"ahler quotient is
\[
  G = S( U(1) \times \ldots U(n-1) \times \prod_{i=1}^{r} U(1) \times
  \ldots \times U(n_i) )
\]
which has rank
\[
  \frac{1}{2} (n^2 - 2 + \sum_{i=1}^{r} n_i^2)
\]
So the real dimension of the implosion is 4 times the rank of $G$, in accordance
with Nakajima's picture. The dimensions and symmetry groups therefore work out correctly--we hope to
further investigate this picture in a future work.

\section{Orthosymplectic examples}

For other classical groups we have to revisit the notion of balance, as well
as the prescription for finding the symmetry group of the dual
(see eg \cite{Hanany-Kalveks}, \cite{Gaiotto-Witten}).
In the case of orthosymplectic quivers (where we use the physics
notation $USp(n) = Sp(n/2))$, there are 2 cases to consider:

\smallskip
(i) that of an orthogonal node labelled by $SO(N)$, with neighbours
$USp(N_j) = Sp(N_j /2)$. The balancing condition is
\[
2N = 2 + \sum N_j
\]
where the sum is taken over all nodes adjacent to the $SO(N)$ one.

\smallskip
(ii) a symplectic mode $USp(N)$ with neighbours $SO(N_j)$. Now the balancing condition is
\[
2N = -2 + \sum N_j
\]

Let us consider the $D_n$ case. The quiver defining the nilpotent
variety is a chain with $2n-2$ gauge nodes $SO(2), USp(2), SO(4), \ldots,
USp(2n-2)$ and then a flavour node $SO(2n)$.
The gauge nodes are all balanced, yielding in the orthosymplectic situation a $SO(2n)$ symmetry in the
dual space. 

For the dual of the implosion, we can mimic the construction in the $A_n$
case, removing the flavour node and replacing it with a bouquet of $n$
$SO(2)$ nodes. This keeps the $USp(2n-2)$ gauge node (and the
preceding gauge nodes) balanced, so we still have
an $SO(2n)$-symmetry in the implosion as required. The unbalanced nodes now
yield a $T^n$ symmetry in the implosion, which again is correct. As in the
$A_n$ case, we have no flavour nodes, reflecting the fact we do not
expect deformation parameters in the implosion.

We can carry out a check using the calculations of Zhenghao Zhong \cite{ZZ}
of the Hilbert series for the Coulomb branch of these quivers for
$n=3,4,5,6,7$. The $t^2$ coefficient, which is expected to give the
dimension of the global symmetry group, is $18,32,50,72,98$ in this cases.
So in each of these cases we obtain

\[
  2n^2 = n + 2n(2n-1)/2  = {\rm rank}\; SO(2n) + {\rm dim}\; SO(2n)
\]
as expected for the complex dimension of the symmetry group of the $SO(2n)$ implosion.

The rank of the  group by which we are performing the hyperk\"ahler
quotient is $n + 2  \sum_{i=1}^{n-1} i = n^2$, and the real dimension
of the $SO(2n)$-implosion is $4n^2$,
in accordance with our expectation.

\smallskip
%For $B_n$, we start with the quiver for the nilpotent variety of
%the Langlands dual $USp(2n)$. This time we have an additional gauge node
%$SO(2n)$ after the $USp(2n-2)$ node, and then a flavour node $USp(2n)$.
%All $2n-1$ gauge nodes are balanced, yielding an $SO(2n+1)$ symmetry of the
%dual.

%We now replace the flavour node by a bouquet of $n$ $USp(2)=Sp(1)$ nodes.
%This keeps the $SO(2n)$ node balanced, so in the dual we get $SO(2n+1)$ symmetry. And we also have $T^n$ symmetry% from the unbalanced nodes. So the dual
%space has the correct symmetries of the implosion.

%The rank of the group by which we quotient is $2 \sum_{i=1}^{n} i =n(n+1)$,
%and the real dimension of the $SO(2n+1)$ implosion is
%$2( n(2n+1) + n) = 4n(n+1)$ as expected.

\medskip
So for  $D_n$ although the original implosion does not appear
to have a quiver description (ie is non-Lagrangian) the dual {\em does}
 arise as a quiver variety.


\begin{thebibliography}{10}


\bibitem{BD} R. Bielawski and A. S. Dancer, \emph{The geometry and topology of toric hyperkahler manifolds}.
  Comm. Anal. Geom. \textbf{8}  727--760 (2000).

\bibitem{BLP} T. Braden,  A. Licata and N. Proudfoot,
   \emph{Gale duality and Koszul duality}. Adv. Math. 225 no. 4, 2002--2049  
   (2010).
   
\bibitem{CHZ} S. Cremonesi, A. Hanany and A. Zaffaroni, \emph{Monopole
    operators and Hilbert series of Coulomb branches of 3d N=4 gauge theories}
  JHEP \textbf{1} (2014) 005
  
\bibitem{Chriss-G:representation} N.~Chriss and V.~Ginzburg,
  \emph{Representation theory and complex geometry}, Birkhauser,
  Boston (1997).

 
\bibitem{DKS} A. Dancer, F. Kirwan and A. Swann, \emph{Implosion for
  hyperk\"ahler manifolds}, Compositio Mathematica \textbf{149}
  592--630 (2013).

\bibitem{DKS-Seshadri} \bysame, \emph{Implosions and hypertoric
  geometry}, J. Ramanujan Math. Soc. \textbf{28A} (special issue in
  honour of Professor C.S. Seshadri's 80th birthday) 81--122 (2013).

\bibitem{DKS-twistor} \bysame, \emph{Twistor spaces for hyperk\"ahler
  implosions}, J. Differential Geometry \textbf{97} (special issue in memory of Friedrich Hirzebruch)  37-77 (2014).

\bibitem{DKS-Arb} A. Dancer, B. Doran, F. Kirwan and A. Swann, 
\emph{Symplectic and hyperk\"ahler implosion},
in Proceedings of the 2013 Arbeitstagung in memory of Friedrich Hirzebruch.
(Eds. W. Ballmann, C. Blohmann, G. Faltings, 
P. Teichner, D. Zagier)
Progress in Mathematics, Birkh\"auser \textbf{319}
 81--103 (2016).

\bibitem{DKR} A. Dancer, F. Kirwan and M. R\"oser, \emph{Hyperk\"ahler
    implosion and Nahm's equations}
  Communications in Mathematical Physics \textbf{342} 251-301 (2016).
  
\bibitem{DK} B.~Doran and F.~Kirwan, \emph{Towards non-reductive
  geometric invariant theory}, Pure Appl. Math. Q. \textbf{3} no. 1,
  part 3, 61--105 (2007).

\bibitem{FHH}  B. Feng, A. Hanany,  Y-H. He,
  \emph{Counting gauge invariants: the plethystic program}.
  J. High Energy Phys. \textbf{3} (2007)

\bibitem{Gaiotto-Witten} D. Gaiotto and E. Witten, \emph{S-duality of boundary conditions in N=4 super Yang-Mills theory}. Adv. Theor. Math. Phys. 13 (2009), 
no. 3, 721--896.

\bibitem{GinzburgRiche} V. Ginzburg and S. Riche, \emph{Differential
  operators on \( G/U \) and the affine Grassmannian},
  J. Math.Jussieu  vol 14    493--575  (2015)
  %preprint
  %\url{arXiv:1306.6754(math.RT)}

\bibitem{GJS:implosion} V.~Guillemin, L.~Jeffrey and
  R.~Sjamaar, \emph{Symplectic implosion}, Transformation Groups
  \textbf{7} 155--184 (2002).

\bibitem{Hanany-Kalveks}
 A. Hanany and R. Kalveks,\emph{Quiver theories for moduli spaces of classical group nilpotent orbits}. J. High Energy Phys. 2016, no. 6, 130, 
front matter+60 pp

\bibitem{K} F. Kirwan, \emph{Symplectic implosion and nonreductive quotients}.
  Geometric aspects of analysis and mechanics, 213-–256, Progr. Math., 292,
  Birkhäuser/Springer, New York, 2011.

\bibitem{Kobak-S:finite} P.~Kobak and A.~Swann, \emph{Classical
  nilpotent orbits as hyperk\"ahler quotients}, International J. Math
  \textbf{7} 193--210 (1996).

\bibitem{Kostant:polynomial} B.~Kostant, \emph{Lie group
  representations on polynomial rings}, Amer. J. Math. \textbf{85}
  327--404 (1963).

\bibitem{KP} H.~Kraft and C.~Procesi, \emph{Minimal singularities in
  \( GL_n \)}, Invent. Math. \textbf{62} 503--515 (1981)

\bibitem{Kronheimer:ALE} P.~B.~Kronheimer,\emph{The construction of ALE spaces as hyper-Kahler quotients},
J. Differential Geom. \textbf{29} 665--683 (1989)
  
\bibitem{Kronheimer:cotangent} P.~B.~Kronheimer, \emph{A hyperk\"ahler
  structure on the cotangent bundle of a compact Lie group}, 1986
  preprint: \url{arXiv:math/0409253} (DG)

\bibitem{MT} G. Moore and Y. Tachikawa, \emph{On 2-D TQFTs whose values are holomorphic symplectic varieties}, in Proceedings of Symposium on Pure Mathematics, String-Math 2011, vol 85 (American Mathematical Society, Providence, RI, 2012)
pp 191-207.

\bibitem{N1} H. Nakajima, \emph{Introduction to a provisional mathematical
    definition of Coulomb branches of 3-dimensional N=4 gauge theories}, in
  Modern geometry: a celebration of the work of Simon Donaldson, 193–-211, Proc. Sympos. Pure Math., 99, Amer. Math. Soc., Providence, RI, (2018).

    
    
\bibitem{GIT} D.~Mumford, J.~Fogarty and F.~Kirwan, \emph{Geometric invariant theory}, 3rd edition, Springer (Berlin, New York) 1994.

\bibitem{Swann} A. Swann, \emph{Hyper-K\"ahler and quaternionic K\"ahler
    geometry}.
  Math. Ann. \textbf{289} 421-450 (1991).
  

\bibitem{ZZ} Z. Zhong \emph{Quiver varieties in 3d, 5d and 6d} MSc thesis, Imperial College, London (2018):
  
  https://imperialcollegelondon.app.box.com/s/gauk9uye1uce8vok5hzi3u5j8jyl7spn

 


\end{thebibliography}
\end{document}